\font\script=eusm10.
\font\sets=msbm10.
\font\stampatello=cmcsc10.
\font\symbols=msam10.

\def\square{\hbox{\vrule\vbox{\hrule\phantom{s}\hrule}\vrule}}
\def\defineq{\buildrel{def}\over{=}}
\def\defin{\buildrel{def}\over{\Longleftrightarrow}}
\def\1{{\bf 1}}

\def\avesum{\sum_{x\sim N}}

\def\C{\hbox{\sets C}}
\def\N{\hbox{\sets N}}

\def\R{\hbox{\sets R}}

\def\Corr{\hbox{\script C}}
\def\EssBdd{\hbox{\symbols n}\,}
\def\modSel{{\widetilde{J}}}

\def\Res{\mathop{{\rm Res}\,}}

\par
\centerline{\bf On the Selberg integral of the three-divisor function $d_3$}
\bigskip
\centerline{\stampatello Giovanni Coppola}
\bigskip
\par
\noindent \centerline{\bf 1. Introduction and statement of the results.}
\smallskip
\par
\noindent
Recall that $d_3(n)\defineq \sum_{abc=n}1$ is the number of ways to write $n$ as a product of three positive integers. Namely, the function $d_3$ is generated by the Dirichlet series $\zeta(s)^3$, where $\zeta$ is the Riemann zeta function.
Given positive integers $N$ and $H=o(N)$ as $N\to \infty$, the {\stampatello Selberg integral} of $d_3$ is
the mean-square 
$$
J_3(N,H)\defineq \sum_{N<x\le 2N} \Big| \sum_{x<n\le x+H}d_3(n)-M_3(x,H)\Big|^2, 
$$
\par
\noindent
where $M_3(x,H)\defineq H\Res_{s=1}\zeta(s)^3x^{s-1}$ is the {\stampatello short intervals mean-value}, i.e. the {\it expected value} of the $d_3-$short sum.  

Improving on results in [CL], we give here a new non-trivial bound for $J_3(N,H)$ by
comparing it to the related {\stampatello modified Selberg integral}, 
$$
\modSel_3(N,H)\defineq \sum_{N<x\le 2N} \Big| \sum_{0\le |n-x|\le H}\Big( 1-{{|n-x|}\over H}\Big)d_3(n)-M_3(x,H)\Big|^2, 
$$
\par
\noindent
where the same choice for the s.i. mean-value $M_3(x,H)$ is due to elementary reasons (see [CL] introduction). 
\par
In what follows, the symbols $O$ and $\ll$ are respectively the usual Landau and Vinogradov notations, while we adopt
the \lq \lq modified Vinogradov notation\rq \rq\ defined as
$$
A(N,H)\EssBdd B(N,H) 
\enspace \defin \enspace 
\forall \varepsilon>0, \enspace A(N,H)\ll_{\varepsilon}N^{\varepsilon} B(N,H). 
$$
\par
\noindent
As any divisor function,  $d_3$ 
is bounded asymptotically by every arbitrarily 
small power of the variable, i.e. it satisfies the definition:
$$
f\enspace \hbox{\stampatello essentially \thinspace bounded} 
\enspace \defin \enspace 
\forall \varepsilon>0 \quad f(n)\ll_{\varepsilon} n^{\varepsilon}\ ,
$$
\par
\noindent
and we shortly write $f\EssBdd 1$. According to the \lq \lq modified Vinogradov notation\rq \rq, we use to write $d_3\EssBdd 1$  when $d_3$ is restricted to $]N-H,2N+H]$ as before. 

\medskip

\par
The author [C] has proved the lower bound $NH\log^4N\ll J_3(N,H)$ for $H\ll N^{1/3}$, while, in an attempt to establish a non trivial upper bound, Laporta and the author
[CL] have conjectured the following estimate for the modified Selberg integral. 
\medskip
\par
\noindent {\stampatello Conjecture CL.} {\it If $H\ll N^{1/3}$, then} $\modSel_3(N,H)\EssBdd NH$.
\medskip
\par
\noindent
As a consequence one has the main result of the present paper. 
\smallskip
\par
\noindent {\stampatello Theorem.} {\it If Conjecture CL holds, then $
J_3(N,H)\EssBdd NH^{6/5}$.}
\medskip
\par
\noindent
The proof of the Theorem is given
in $\S3$, where Conjecture CL is combined with the following general Proposition on
the Selberg integral and the modified one of an essentially bounded arithmetic function $f$, respectively
$$
J_f(N,H)\defineq \avesum \Big| \sum_{x<n\le x+H}f(n) - M_f(x,H)\Big|^2, 
$$
$$
\modSel_f(N,H)\defineq \avesum \Big| \sum_{0\le |n-x|\le H}\Big( 1-{{|n-x|}\over H}\Big)f(n) - M_f(x,H)\Big|^2\ .
$$
\par
\noindent
Here $x\sim N$ means that $N<x\le 2N$ and the s.i. mean-value is defined as
$$
M_f(x,H)\defineq H\Res_{s=1}F(s)x^{s-1},
$$
\par
\noindent
whenever $f$ is generated by a meromorphic Dirichlet series $F$ with (at most) a pole in $s=1$.
It is not difficult to see that $M_f(x,H)\sim \sum_{x<n\le x+H}p_c(\log n)$ with good remainders, 
where $p_c$ is the so called {\stampatello logarithmic polynomial} of $f$ of degree $c={\rm ord}_{s=1}F(s)-1$ (see [CL]).
If $M_f(x,H)$ vanishes identically (that is the case when $F$ is regular at $s=1$), then we say that $f$ is {\stampatello balanced}.
\medskip
\par
\noindent 
{\stampatello Proposition.} {\it If
$f:\N \rightarrow \C$ is an essentially bounded and balanced function such that 
$$
\modSel_f(N,H)\EssBdd NH^{1+A}, \quad \forall H\in [H_1,H_2]\
\hbox{with}\ N^{\delta}\ll H_1\ll H_2\ll N^{1/2-\delta}\ ,
$$
\par
\noindent 
for an absolute constant $A\in [0,1)$ and 
for a fixed $\delta\in(0,1/2)$,
then} 
$$
J_f(N,H)\EssBdd NH^{1+{{1+3A}\over {5-A}}}, \quad \forall H\le H_2\ . 
$$

Beyond Lemma 1 of [CL] (whose formul\ae\ are quoted in \S2), for the proof of the Proposition (see $\S3$) it is applied an enhanced version of Gallagher's Lemma for the exponential sums (see [Ga], Lemma 1) established by the author in the joint paper [CL1] with Laporta. 
Unfortunately at the moment Conjecture CL  remains unproved after a serious gap occurred in the proof of the \lq \lq Fundamental Lemma\rq \rq \thinspace given in a former version of [CL]. Nevertheless, we think it is worthwhile to make further attempts, 
first for the strongest conjecture \thinspace $J_3(N,H)\EssBdd NH$, $\forall H\ll N^{1/3}$ seems to be out of reach by any method, then because the bound of the Theorem has an impressive consequence on the $6$-th moment of $\zeta$ on the critical line (compare [CL] and version 2 of the present paper on arXiv). 
Actually, the author has already explored in his paper [C2]
a link between the Selberg integral $J_k(N,H)$ of the $k$-divisor function $d_k$
and the $2k-$th moment of Riemann zeta function on the critical line, i.e. 
$$
I_k(T)\defineq \int_{T}^{2T}\Big| \zeta\Big( {1\over 2}+it\Big)\Big|^{2k} {\rm d}t 
$$
\par
\noindent
(under suitably conditions on $N,H,T$). Somehow the same link holds also
between $\modSel_k(N,H)$ and $I_k(T)$, since in Gallagher's Lemma the right hand side
is a Selberg integral de facto and, as showed in [CL1], it
can be replaced by the corresponding
modified Selberg integral (see Lemma in \S2). This is still true for the
Dirichlet polynomials case, namely Theorem 1 of [Ga] can be modified in the same fashion. 
Thus, by applying such a Dirichlet polynomial version of Gallagher's Lemma within the method of [C2] instead of Theorem 1 of [Ga] one gets the following outstanding consequence of Conjecture CL.
\smallskip
\par
\noindent 
{\stampatello Corollary.} {\it If Conjecture CL holds, then} $
I_3(T)\EssBdd T$.
\medskip
\par
\noindent 
In the literature this result is known as the \lq \lq weak $6-$th moment\rq \rq, since it gives no asymptotic equality for $I_3(T)$, but just  an upper bound with additional sufficiently small powers. The proof of the Corollary will be given in a forthcoming paper. 
\bigskip

\par
\noindent \centerline{\bf 2. Notation and preliminary formul\ae.}
\smallskip
\par
\noindent
The {\stampatello correlation} of $f$ with {\stampatello shift} $h$ is defined as 
$$
\Corr_f(h)\defineq \sum_{n\sim N}f(n)\overline{f(n-h)}=\sum_{n\sim N}f(n)\sum_{{m\sim N}\atop {n-m=h}}\overline{f(m)}+O\Big( \max_{[N-|h|,2N+|h|]}|f|^2\sum_{n\in [N-|h|,N]\cup [2N,2N+|h|]}1\Big).
$$
\par
\noindent
Recall that from the orthogonality of the exponentials $e(\beta)\defineq e^{2\pi i\beta}$ ($\beta \in \R$) one has
$$
\sum_{n\sim N}f(n)\sum_{{m\sim N}\atop {n-m=h}}\overline{f(m)}=\int_{0}^{1}| \widehat{f}(\alpha)|^2 e(-h\alpha){\rm d}\alpha, 
$$
\par				
\noindent
where
$$
\widehat{f}(\alpha)\defineq \sum_n f(n)e(n\alpha) = \sum_{n\sim N} f(n)e(n\alpha)
$$
is truncated in the range $]N,2N]$ in order to avoid convergence problems. 
If $f\EssBdd 1$,  then (compare [CL])
$$
\Corr_f(h)=\int_{0}^{1}| \widehat{f}(\alpha)|^2 e(-h\alpha){\rm d}\alpha+O_{\varepsilon}\Big( N^{\varepsilon}|h|\Big),\quad  
\forall h\not=0.
$$
\par
\noindent

The {\stampatello correlation} of an uniformly bounded {\it weight} $w:\R\rightarrow \C$, vanishing outside $[-H,H]$, is conveniently defined as 
$$
\Corr_w(h)\defineq \sum_{a}w(a)\sum_{b\atop {a-b=h}}\overline{w(b)}. 
$$
\par
\noindent
In particular, the main weight involved here is the characteristic function of the integers in $[1,H]$, say $u(n)\defineq \1_{[1,H]}(n)$, whose  correlation is
$$
\Corr_u(h)\defineq \sum_{a}u(a)\sum_{b\atop {a-b=h}}u(b) = \sum_{a\le H}\sum_{{b\le H}\atop {a-b=h}}1 = \max(H-|h|,0). 
$$
\par
\noindent
If $f$ is essentially bounded and balanced, then
Lemma 1 of [CL] provides the formul\ae\ 
$$
J_f(N,H)=\sum_h \Corr_u(h)\Corr_f(h)+O_{\varepsilon}( N^{\varepsilon}H^3),
\qquad
\modSel_f(N,H)=\sum_h \Corr_{\Corr_u/H}(h)\Corr_f(h)+O_{\varepsilon}( N^{\varepsilon}H^3),
$$
\par
\noindent
where we find the {\it Cesaro weight} ${\displaystyle  {\Corr_u(h)\over H}=\max\Big(1-{{|h|}\over H},0\Big)}$, that is the \lq \lq normalized\rq \rq \thinspace correlation of $u$.\par\noindent
In what follows such formul\ae\ will be applied in the integral form
$$
\modSel_f(N,H)=\int_{-1/2}^{1/2}\left|\widehat{f}(\alpha)\right|^2 {{|\widehat{u}(\alpha)|^4}\over {H^2}} {\rm d}\alpha 
 + O_{\varepsilon}(N^{\varepsilon}H^3), 
\quad 
J_f(N,H)=\int_{-1/2}^{1/2}\left|\widehat{f}(\alpha)\right|^2 |\widehat{u}(\alpha)|^2 {\rm d}\alpha 
 + O_{\varepsilon}(N^{\varepsilon}H^3)\ ,
$$
obtained by using the exponential sum $\widehat{u}(\alpha)\defineq \sum_{h\le H}e(h\alpha)$.
Note that for every $0<\varepsilon<1$ one has
$$
|\widehat{u}(\alpha)|
>[\varepsilon H]
\enspace \Longrightarrow \enspace
|\alpha|< {1\over {2[\varepsilon H]}}. 
\leqno{(\ast)}
$$
\par\noindent
Indeed, since $(\ast)$ is trivial for $\alpha=0$, we may
assume that $0<\alpha<1/2$, which easily implies
$2\alpha<\sin(\pi\alpha)$. Then, $(\ast)$ follows also in this case for
$$
[\varepsilon H]<|\widehat{u}(\alpha)|
={{|\sin(\pi H \alpha)|}\over {\sin(\pi \alpha)}}\le
{{1}\over { \sin(\pi \alpha)}}< {{1}\over {2\alpha}}.
$$
Finally, we quote the aforementioned modified version of Gallagher's Lemma (see [CL1]).
\smallskip
\par
\noindent {\stampatello Lemma.}
{\it Let \thinspace $N,h$ \thinspace be positive integers such that \thinspace $h\to \infty$ \thinspace and \thinspace $h=o(N)$ as \thinspace $N\to \infty$. If \thinspace $f:\N \rightarrow \C$ \thinspace is essentially bounded and balanced,
then}
$$
h^2\int_{-{1\over {2h}}}^{{1\over {2h}}}\Big| \widehat{f}(\alpha)\Big|^2 {\rm d}\alpha
\EssBdd \modSel_f(N,h)+h^3. 
$$

\bigskip

\par
\noindent \centerline{\bf 3. Proofs of the Proposition and the Theorem.}
\smallskip
\par

\noindent {\stampatello Proof of the Proposition.} If $H$ were not essentially bounded, the trivial bound $\EssBdd NH^2$ would imply $J_f(N,H)\EssBdd N$ immediately. Thus, let us assume that $H\gg N^{\eta}$ for some small $\eta>0$.
Taking $\varepsilon=\varepsilon(H), E=E(H)$ to be determined later such that
$0<\varepsilon<E$ and $\varepsilon, E \to 0$ as $H\to\infty$, we write
$$
J_f(N,H)\EssBdd \int_{-1/2}^{1/2}\left|\widehat{f}(\alpha)\right|^2 |\widehat{u}(\alpha)|^2 {\rm d}\alpha + H^3\EssBdd
$$
$$
\EssBdd \varepsilon^2 H^2 \int_{|\widehat{u}(\alpha)|\le [\varepsilon H]}\left|\widehat{f}(\alpha)\right|^2 {\rm d}\alpha
 + E^2 H^2 \int_{[\varepsilon H]<|\widehat{u}(\alpha)|\le EH}\left|\widehat{f}(\alpha)\right|^2 {\rm d}\alpha
 + {1\over {E^2}} \int_{|\widehat{u}(\alpha)|>EH}\left|\widehat{f}(\alpha)\right|^2 {{|\widehat{u}(\alpha)|^4}\over {H^2}} {\rm d}\alpha + H^3.
$$
\par
\noindent
By applying Parseval's identity together with $(\ast)$ one has
$$
J_f(N,H)\EssBdd NH^2 \varepsilon^2 + H^2 E^2 \int_{|\alpha|\le {1\over {2h}}}\left|\widehat{f}(\alpha)\right|^2 {\rm d}\alpha
 + {1\over {E^2}} \int_{-1/2}^{1/2}\left|\widehat{f}(\alpha)\right|^2 {{|\widehat{u}(\alpha)|^4}\over {H^2}} {\rm d}\alpha + H^3\EssBdd
$$
$$
\EssBdd NH^2 \varepsilon^2 + H^2 E^2 \int_{|\alpha|\le {1\over {2h}}}\left|\widehat{f}(\alpha)\right|^2 {\rm d}\alpha 
 + {1\over {E^2}}\modSel_f(N,H) + {1\over {E^2}}H^3. 
$$
\par
\noindent
Now the hypothesis \enspace $\modSel_f(N,h)\EssBdd Nh^{1+A}$ \enspace and the previous Lemma for $h\defineq [\varepsilon H] \to \infty$ imply
$$
J_f(N,H)\EssBdd NH^2 \Big( \varepsilon^2 + {{E^2 H^A}\over {\varepsilon^{1-A} H}} + {{H^A}\over {E^2 H}}\Big). 
$$
\par
\noindent
Taking \thinspace $\varepsilon = H^{-{{2(1-A)}\over {5-A}}}$, \thinspace $E = H^{-{{(1-A)^2}\over {2(5-A)}}}$ one gets
\thinspace $\varepsilon^2 = {{E^2 H^{A-1}}\over {\varepsilon^{1-A}}} = {{H^{A-1}}\over {E^2}}$, $\varepsilon=o(E)$ and $E\to 0$, as desired.\enspace $\square$ 

\medskip

\par
\noindent {\stampatello Proof of the Theorem.} 
Denoting the logarithmic polynomial of $d_3$ with $p_2$, we recall that $M_3(x,H)\sim \sum_{x<n\le x+H}p_2(\log n)$ where the implicit remainders give a negligible contribution to $J_3(N,H)$. Thus,
it is sufficient to apply the Proposition to the balanced and essentially bounded function $f(n)=d_3(n)-p_2(\log n)$, because Conjecture CL allows to take $A=0$ and $H_2\ll N^{1/3}$.\enspace $\square$ 
\bigskip

\par
\noindent {\bf Acknowledgement}. The author wishes to thank Maurizio Laporta for useful comments. 

\bigskip

\par
\centerline{\stampatello References}
\medskip
\item{\bf [C]} Coppola, G. \thinspace - \thinspace {\sl On some lower bounds of some symmetry integrals} \thinspace - \thinspace http://arxiv.org/abs/1003.4553 \thinspace - \thinspace to appear on Afrika Mathematika (Springer) 
\item{\bf [C2]} Coppola, G. \thinspace - \thinspace {\sl On the Selberg integral of the $k$-divisor function and the $2k$-th moment of the Riemann zeta-function} Publ. Inst. Math. (Beograd) (N.S.) {\bf 88(102)} (2010), 99--110. \thinspace - \thinspace available online
\item{\bf [CL]} Coppola, G. and Laporta, M. \thinspace - \thinspace {\sl Generations of correlation averages} \thinspace - \thinspace http://arxiv.org/abs/1205.1706 
\item{\bf [CL1]} Coppola, G. and Laporta, M. \thinspace - \thinspace {\sl A modified Gallagher's Lemma} \thinspace - \thinspace http://arxiv.org/abs/ 
\item{\bf [Ga]} Gallagher, P. X. \thinspace - \thinspace {\sl A large sieve density estimate near $\sigma =1$} \thinspace - \thinspace Invent. Math. {\bf 11} (1970), 329--339. $\underline{\tt MR\enspace 43\# 4775}$

\bigskip

\par

\leftline{\tt Giovanni Coppola}
\leftline{\tt Universit\`a degli Studi di Salerno}
\leftline{\tt Home address \negthinspace : \negthinspace Via Partenio (Pal.Soldatiello) - 83100, Avellino(AV) - ITALY}
\leftline{\tt e-page : $\! \! \! \! \! \!$ www.giovannicoppola.name}
\leftline{\tt e-mail : gcoppola@diima.unisa.it}

\bye